\documentclass[11pt]{article}
\usepackage{amssymb,amsmath}

\newcommand{\be}{\begin{equation}}
\newcommand{\ee}{\end{equation}}
\newcommand{\bea}{\begin{eqnarray}}
\newcommand{\eea}{\end{eqnarray}}
\newtheorem{theorem}{Theorem}

\newtheorem{corollary}{Corollary}
\newtheorem{example}{Example}
\newtheorem{lemma}{Lemma}

\def\1#1{^{(#1)}}

\begin{document}
\title{Symmetric Numerical Semigroups \\Generated by Fibonacci and Lucas 
Triples}
\author{Leonid G. Fel\\
\\Department of Civil Engineering, Technion, Haifa 3200, Israel\\
\vspace{-.2cm}
\\{\sl e-mail: lfel@tx.technion.ac.il}}
\date{\today}
\maketitle
\def\be{\begin{equation}}
\def\ee{\end{equation}}
\def\bea{\begin{eqnarray}}
\def\eea{\end{eqnarray}}
\def\p{\prime}
\vspace{-.2cm}
\begin{abstract}
The symmetric numerical semigroups ${\sf S}\left(F_a,F_b,F_c\right)$ and ${\sf 
S}\left(L_k,L_m,L_n\right)$ generated by three Fibonacci $(F_a,F_b,F_c)$ and  
Lucas $(L_k,L_m,L_n)$ numbers are considered. Based on divisibility properties 
of the Fibonacci and Lucas numbers we establish necessary and sufficient 
conditions for both semigroups to be symmetric and calculate their Hilbert 
generating series, Frobenius numbers and genera.
\end{abstract}
\section{Introduction}\label{s1}
Recently the numerical semigroups ${\sf S}\left(F_i,F_{i+2},F_{i+k}\right)$, 
$i,k\geq 3$, generated by three Fibonacci numbers $F_j$ were discussed in 
\cite{mrr07}. It turns out that the remarkable properties of $F_j$ in these 
triples suffice to calculate the Frobenius number ${\cal F}\left({\sf S}\right)$
and genus $G\left({\sf S}\right)$ of semigroup. In this article we show that a 
nature of Fibonacci and Lucas numbers is sufficient not only to calculate the 
specific parameters of semigroups, but also to describe completely the structure
of symmetric numerical semigroups ${\sf S}\left(F_a,F_b,F_c\right)$, $3\leq a<b
<c$, and ${\sf S}\left(L_k,L_m,L_n\right)$, $2\leq k<m<n$, generated by 
Fibonacci
\footnote{We avoid to use the term "{\em Fibonacci semigroup}" because it has 
been already reserved for another algebraic structure \cite{res89}.}
and Lucas numbers, respectively. Based on divisibility properties of these 
numbers we establish necessary and sufficient conditions for both semigroups to 
be symmetric and calculate their Hilbert generating series, Frobenius numbers 
and genera.
\section{Basic properties of the 3D symmetric numerical semigroups}\label{s2}
Recall basic definitions and known facts about 3D numerical semigroups mostly 
focusing on their symmetric type. Let ${\sf S}\left(d_1,d_2,d_3\right)\subset 
{\mathbb Z}_+\cup\{0\}$ be the additive numerical semigroup with zero finitely 
generated by a minimal set of positive integers $\{d_1,d_2,d_3\}$ such that $3
\leq d_1<d_2<d_3$, $\gcd(d_1,d_2,d_3)=1$. Semigroup ${\sf S}(d_1,d_2,d_3)$ is 
said to be generated by the minimal set of three natural numbers if there are 
no nonnegative integers $b_{i,j}$ for which the following dependence holds:
\begin{equation}
d_i=\sum_{j\neq i}^mb_{i,j}d_j\;,\;\;\;b_{i,j}\in\{0,1,\ldots\}\;\;\; 
\mbox{for any}\;\;i\leq m\;.\label{deff1}
\end{equation}
For short we denote the vector $(d_1,d_2,d_3)$ by ${\bf d}^3$. Following Johnson
\cite{john60} define {\em the minimal relation} ${\cal R}_3$ for given ${\bf 
d}^3$ as follows
\begin{eqnarray}
{\cal R}_3\left(\begin{array}{r}d_1\\d_2\\d_3 \end{array}\right)=
\left(\begin{array}{r}0\\0\\0\end{array}\right)\;,\;\;\;
{\cal R}_3=\left(\begin{array}{rrr}a_{11} & -a_{12} & -a_{13} \\
-a_{21} & a_{22} &-a_{23}\\-a_{31}& -a_{32}&a_{33}\end{array}\right)\;,\;\;\;
\left\{\begin{array}{r}\gcd(a_{11},a_{12},a_{13})=1\\
\gcd(a_{21},a_{22},a_{23})=1\\\gcd(a_{31},a_{32},a_{33})=1\end{array}\right.\;,
\label{herznon2}
\end{eqnarray}
where
\begin{eqnarray}
a_{11}&=&\min\left\{v_{11}\;|\;v_{11}\geq 2,\;v_{11}d_1=v_{12}d_2+v_{13}d_3,\;
v_{12},v_{13}\in {\mathbb N}\cup\{0\}\right\}\;,\nonumber\\
a_{22}&=&\min\left\{v_{22}\;|\;v_{22}\geq 2,\;v_{22}d_2=v_{21}d_1+v_{23}d_3,\;
v_{21},v_{23}\in {\mathbb N}\cup\{0\}\right\}\;,\label{herznon3}\\
a_{33}&=&\min\left\{v_{33}\;|\;v_{33}\geq 2,\;v_{33}d_3=v_{31}d_1+
v_{32}d_2,\;v_{31},v_{32}\in {\mathbb N}\cup\{0\}\right\}\;.\nonumber
\end{eqnarray}
The uniquely defined values of $v_{ij},i\neq j$ which give $a_{ii}$ will be
denoted by $a_{ij},i\neq j$. Note that due to minimality of the set $(d_1,
d_2,d_3)$ the elements $a_{ij},i,j\leq 3$ satisfy
\begin{eqnarray}
&&a_{11}=a_{21}+a_{31}\;,\;\;\;a_{22}=a_{12}+a_{32}\;,\;\;\;
a_{33}=a_{13}+a_{23}\;,\nonumber\\
&&d_1=a_{22}a_{33}-a_{23}a_{32}\;,\;\;\;d_2=a_{11}a_{33}-a_{13}a_{31}\;,\;\;\;
d_3=a_{11}a_{22}-a_{12}a_{21}\;.\label{herz4}
\end{eqnarray}
The smallest integer $C\left({\bf d}^3\right)$ such that all integers $s,\;s
\geq C\left({\bf d}^3\right)$, belong to ${\sf S}\left({\bf d}^3\right)$ is 
called {\em the conductor} of ${\sf S}\left({\bf d}^3\right)$,
\begin{eqnarray}
C\left({\bf d}^3\right):=\min\left\{s\in {\sf S}\left({\bf d}^3\right)\;|\;s+
{\mathbb Z}_+\cup\{0\}\subset {\sf S}\left({\bf d}^3\right)\right\}\;.\nonumber
\end{eqnarray}
The number ${\cal F}\left({\bf d}^3\right)=C\left({\bf d}^3\right)-1$ is 
referred to as {\em the Frobenius number}. Denote by $\Delta\left({\bf d}^3
\right)$ the complement of ${\sf S}\left({\bf d}^3\right)$ in ${\mathbb Z}_+
\cup\{0\}$, i.e. $\Delta\left({\bf d}^3\right)={\mathbb Z}_+\cup\{0\}\setminus 
{\sf S}\left({\bf d}^3\right)$. The cardinality ($\#$) of the set $\Delta
\left({\bf d}^3\right)$ is called {\em the number of gaps}, $G\left({\bf d}^3
\right):=\#\left\{\Delta\left({\bf d}^3\right)\right\}$, or {\em genus} of 
${\sf S}\left({\bf d}^3\right)$. 

The semigroup ring ${\sf k}\left[X_1,X_2,X_3\right]$ over a field ${\sf k}$ of 
characteristic 0 associated with ${\sf S}\left({\bf d}^3\right)$ is a polynomial
subring graded by $\deg X_i=d_i$, $i=1,2,3$ and generated by all monomials 
$z^{d_i}$. The Hilbert series $H({\bf d}^3;z)$ of a graded subring ${\sf k}
\left[z^{d_1},z^{d_2},z^{d_3}\right]$ is defined \cite{stan96} by
\begin{equation}
H({\bf d}^3;z)=\sum_{s\;\in\;{\sf S}\left({\bf d}^3\right)} z^s=\frac{Q({\bf 
d}^3;z)}{\left(1-z^{d_1}\right)\left(1-z^{d_2}\right)\left(1-z^{d_3}\right)}\;,
\label{hilb0}
\end{equation}
where $Q({\bf d}^3;z)$ is a polynomial in $z$.

The semigroup ${\sf S}\left({\bf d}^3\right)$ is called {\em symmetric} iff  
for any integer $s$ holds
\begin{eqnarray}
s\in {\sf S}\left({\bf d}^3\right)\;\;\;\Longleftrightarrow\;\;\;{\cal F}\left(
{\bf d}^3\right)-s\not\in{\sf S}\left({\bf d}^3\right)\;.\label{intro3}
\end{eqnarray}
Otherwise ${\sf S}\left({\bf d}^3\right)$ is called {\em non--symmetric}. The
integers $G\left({\bf d}^3\right)$ and $C\left({\bf d}^3\right)$ are related
\cite{heku71} as,
\begin{eqnarray}
2G\left({\bf d}^3\right)=C\left({\bf d}^3\right)\;\;\mbox{if}\;\;{\sf S}
\left({\bf d}^3\right)\;\;\mbox{is symmetric semigroup, and}\;\;2G\left({\bf
d}^3\right)>C\left({\bf d}^3\right)\;\;\mbox{otherwise}.\label{intro5}
\end{eqnarray}
Notice that ${\sf S}\left({\bf d}^2\right)$ is always symmetric semigroup
\cite{aper46}. The number of independent entries $a_{ij}$ in (\ref{herznon2}) 
can be reduced if ${\sf S}\left({\bf d}^3\right)$ is symmetric: at least one 
off-diagonal element of $\widehat {\cal R}_3$ vanishes, e.g. $a_{13}=0$ and 
therefore $a_{11}d_1=a_{12}d_2$. Due to {\em minimality} of the last relation 
we have by (\ref{herznon2}) the following equalities and consequently the matrix
representation as well \cite{herz70} (see also \cite{fel04}, Section 6.2)
\begin{eqnarray}
\left.\begin{array}{l}a_{11}=a_{21}={\sf lcm}(d_1,d_2)/d_1,\;\;\;
a_{12}=a_{22}={\sf lcm}(d_1,d_2)/d_2\;,\\
a_{33}=d_1/a_{22}=d_2/a_{11}\;,\;\;a_{23}=0\;,\end{array}\right.\;\;
\widehat {\cal R}_{3s}=\left(\begin{array}{rrr}a_{11} & -a_{22} & 0 \\
-a_{11} & a_{22} & 0 \\-a_{31} & -a_{32} & a_{33} \end{array}\right),
\label{herz5}
\end{eqnarray}
where subscript "$s$" stands for symmetric semigroup. Combining (\ref{herz5}) 
with formula for the Frobenius number of symmetric semigroup \cite{herz70}, 
${\cal F}\left({\bf d}^3_s\right)=a_{22}d_2+a_{33}d_3-\sum_{i=1}^3d_i$, we get 
finally, 
\begin{eqnarray}
{\cal F}\left({\bf d}^3_s\right)=e_1+e_2-\sum_{i=1}^3d_i\;,\;\;\;\;e_1={\sf lcm}
(d_1,d_2)\;,\;\;\;e_2=d_3\; {\sf gcd}(d_1,d_2)\;.\label{herz6}
\end{eqnarray}

If ${\sf S}\left({\bf d}^3\right)$ is symmetric semigroup then ${\sf k}\left[
{\sf S}\left({\bf d}^3\right)\right]$ is a complete intersection \cite{herz70} 
and the numerator $Q({\bf d}^3;z)$ in the Hilbert series (\ref{hilb0}) reads 
\cite{stan96}
\begin{eqnarray}
Q({\bf d}^3;z)=(1-z^{e_1})(1-z^{e_2})\;.\label{sylves10}
\end{eqnarray}
\subsection{Structure of generating triples of symmetric numerical semigroups}
\label{s2a}
Two following statements, Theorem \ref{the1} and Corollary \ref{cor1}, give 
necessary and sufficient conditions for ${\sf S}\left({\bf d}^3\right)$ to be 
symmetric.
\begin{theorem}{\rm (\cite{herz70} and Proposition 3, \cite{wata73})}
\label{the1}
If a semigroup ${\sf S}\left(d_1,d_2,d_3\right)$ is symmetric then its minimal
generating set has the following presentation with two relatively not prime 
elements:
\begin{eqnarray}
\gcd(d_1,d_2)=\lambda\;,\;\;\gcd(d_3,\lambda)=1\;,\;\;d_3\in {\sf S}\left(
\frac{d_1}{\lambda},\frac{d_2}{\lambda}\right)\;.\label{wat2}
\end{eqnarray}
\end{theorem}
It turns out that (\ref{wat2}) gives also sufficient conditions for ${\sf S}
\left({\bf d}^3\right)$ to be symmetric. This follows by Corollary \ref{cor1} of
the old Lemma of Watanabe \cite{wata73} for semigroup ${\sf S}\left({\bf d}^m
\right)$
\begin{lemma}{\rm (Lemma 1, \cite{wata73})}\label{lem1}
Let ${\sf S}\left(d_1,\ldots,d_m\right)$ be a numerical semigroup, $a$ and $b$
be positive integers such that: (i) $c\in{\sf S}\left(d_1,\ldots,d_m\right)$
and $c\neq d_i$, (ii) $\gcd(c,\lambda)=1$.\\
Then semigroup ${\sf S}\left(\lambda d_1,\ldots,\lambda d_m,c\right)$ is 
symmetric iff ${\sf S}\left(d_1,\ldots,d_m\right)$ is symmetric.
\end{lemma}
Combining Lemma \ref{lem1} with the fact that every semigroup ${\sf S}
\left({\bf d}^2\right)$ is symmetric we arrive at Corollary.
\begin{corollary}\label{cor1}
Let ${\sf S}\left(d_1,d_2\right)$ be a numerical semigroup, $c$ and $\lambda$ be
positive integers, $\gcd(c,\lambda)=1$. If $c\in{\sf S}\left(d_1,d_2\right)$, 
then the semigroup ${\sf S}\left(\lambda d_1,\lambda d_2,c\right)$ is symmetric.
\end{corollary}
In Corollary \ref{cor1} the requirement $c\neq d_1,d_2$ can be omitted since
both semigroups ${\sf S}\left(\lambda d_1,\lambda d_2,d_1\right)$ and 
${\sf S}\left(\lambda d_1,\lambda d_2,d_2\right)$ are generated by two 
elements ($d_1,\lambda d_2$) and are also symmetric.

Finish this Section with important proposition adapted to the 3D numerical 
semigroups.
\begin{theorem}{\rm (\cite{heku71}, Proposition 1.14)}\label{the2}\\
The numerical semigroup ${\sf S}\left(3,d_2,d_3\right)$, $\gcd(3,d_2,d_3)=1$, 
$3\nmid d_2$ and $d_3\not\in{\sf S}\left(3,d_2\right)$, is never symmetric.
\end{theorem}
\section{Divisibility of Fibonacci and Lucas numbers}\label{s3}
We recall a remarkable divisibility properties of Fibonacci and Lucas numbers 
which are necessary for further consideration. Theorem \ref{the3} dates back to
E. Lucas \cite{luc78} (Section 11, p. 206),
\begin{theorem}\label{the3}
Let $F_m$ and $F_n$, $m>n$, be the Fibonacci numbers. Then 
\begin{eqnarray}
\gcd\left(F_m,F_n\right)=F_{\gcd(m,n)}\;.\label{wat3}
\end{eqnarray} 
\end{theorem}

As for Theorem \ref{the4}, its weak version was given by Carmichael 
\cite{car13}  
\footnote{Carmichael \cite{car13} (Theorem 7, p. 40) has proven only the most 
hard part of Theorem \ref{the4}, namely, the 1st equality in (\ref{wat3a}).}.
We present here its modern form proved by 
Ribenboim \cite{rib89} and McDaniel \cite{wmd91}.
\begin{theorem}\label{the4}
Let $L_m$ and $L_n$ be the Lucas numbers, and let $m=2^am^{\p}$, $n=2^bn^{\p}$,
where $m^{\p}$ and $n^{\p}$ are odd positive integers and $a,b\geq 0$. Then
\begin{eqnarray}
\gcd\left(L_m,L_n\right)=\left\{\begin{array}{lll}L_{\gcd(m,n)}&\mbox{if}&
a=b\;,\\2 &\mbox{if} &a\neq b\;,\;\;3\mid \gcd(m,n)\;,\\
1 &\mbox{if} &a\neq b\;,\;\;3\nmid \gcd(m,n)\;.\end{array}\right.\label{wat3a}
\end{eqnarray}
\end{theorem}
We also recall another basic divisibility property of Lucas numbers,
\begin{eqnarray}
L_m=0\pmod 2\;,\;\;\;\;\mbox{iff}\;\;\;\;m=0\pmod 3\;.\label{wat3c}  
\end{eqnarray}
We'll need a technical Corollary which follows by consequence of Theorem 
\ref{the4}.
\begin{corollary}\label{cor2}
Let $L_m$ and $L_n$ be the Lucas numbers, and let $m=2^am^{\p}$, $n=2^bn^{\p}$,
where $m^{\p}$ and $n^{\p}$ are odd positive integers and $a,b\geq 0$. Then
\begin{eqnarray}
\gcd\left(L_m,L_n\right)=1\;,\;\;\;\;\mbox{iff}\;\;\;\;\left\{\begin{array}{ll}
a=b=0\;,&\gcd\left(m^{\p},n^{\p}\right)=1\;,\\
a\neq b\;,&\gcd\left(3,\gcd(m,n)\right)=1\;.\end{array}\right.\label{wat3b}
\end{eqnarray}
\end{corollary}
\section{Symmetric numerical semigroups generated by Fibonacci triple}\label{s4}
In this Section we consider symmetric numerical semigroups generated by three 
Fibonacci numbers $F_c$, $F_b$ and $F_a$, $c>b>a\geq 3$. The two first values 
$a=3,4$ are of special interest because of Fibonacci numbers $F_3=2$ and $F_4=
3$. First, the semigroup ${\sf S}\left(F_3,F_b,F_c\right)$, $\gcd(2,F_b,F_c)=1$,
is always symmetric and has actually 2 generators. Next, according to Theorem 
\ref{the2} the semigroup ${\sf S}\left(F_4,F_b,F_c\right)$ is symmetric iff at 
least one of two requirements, $3\nmid F_b$ and $F_c\not\in{\sf S}\left(3,F_b
\right)$, is broken. Avoiding those trivial cases we state 
\begin{theorem}\label{the5}
Let $F_c$, $F_b$ and $F_a$ be the Fibonacci numbers where $c>b>a\geq 5$. Then a 
numerical semigroup ${\sf S}\left(F_a,F_b,F_c\right)$ is symmetric iff
\begin{eqnarray}
\lambda=\gcd(a,b)\geq 3\;,\;\;\;\gcd(\lambda,c)=1,2\;,\;\;\;F_c\in {\sf S}
\left(\frac{F_a}{F_{\lambda}},\frac{F_b}{F_{\lambda}}\right)\;,\label{wat4}
\end{eqnarray}
\end{theorem}
{\sf Proof} $\;\;\;$By Theorem \ref{the1} and Corollary \ref{cor1} a numerical 
semigroup ${\sf S}\left(F_a,F_b,F_c\right)$ is symmetric iff 
\begin{eqnarray}
g=\gcd\left(F_a,F_b\right)>1\;,\;\;\;\gcd(g,F_c)=1\;,\;\;\;F_c\in {\sf 
S}\left(\frac{F_a}{g},\frac{F_b}{g}\right)\;.\label{wat5}
\end{eqnarray}
By consequence of Theorem \ref{the3} and definition of Fibonacci numbers we get
\begin{eqnarray}
\left\{\begin{array}{lll}g=F_{\lambda}>1 &\rightarrow &\gcd(a,b)\geq 3\;,\\
\gcd(F_{\lambda},F_c)=F_{\gcd(\lambda,c)}=1 &\rightarrow &\gcd(\lambda,c)=1,2
\;.\end{array}\right.\label{wat6}
\end{eqnarray}
The last containment in (\ref{wat5}) gives
\begin{eqnarray}
F_c=A\frac{F_a}{g}+B\frac{F_b}{g}=A\frac{F_a}{F_{\lambda}}+B\frac{F_b}{F_{
\lambda}}\;,\;\;\;\;A,B\in {\mathbb Z}_+\;,\nonumber
\end{eqnarray}
that finishes the proof of Theorem.$\;\;\;\;\;\;\Box$

Theorem \ref{the5} remains true for any permutation of indices in triple $(F_a,
F_b,F_c)$. By (\ref{herz6}), (\ref{sylves10}) and (\ref{wat4}) we get
\begin{corollary}\label{cor3}
Let $F_c$, $F_b$ and $F_a$ be the Fibonacci numbers and numerical 
semigroup ${\sf S}\left(F_a,F_b,F_c\right)$ be symmetric. Then its Hilbert 
series and Frobenius number are given by
\begin{eqnarray}
H\left(F_a,F_b,F_c\right)&=&\frac{(1-z^{f_1})(1-z^{f_2})}{\left(1-z^{F_a}\right)
\left(1-z^{F_b}\right)\left(1-z^{F_c}\right)}\;,\;\;\;f_1=\frac{F_aF_b}{F_{
\gcd(a,b)}}\;,\;\;\;f_2=F_c\cdot F_{\gcd(a,b)}\;,\label{wat8}\\
{\cal F}\left(F_a,F_b,F_c\right)&=&f_1+f_2-(F_a+F_b+F_c)\;.\nonumber
\end{eqnarray}
\end{corollary}
The next Corollary \ref{cor4} gives only the sufficient condition for ${\sf S}
\left(F_a,F_b,F_c\right)$ to be symmetric and is less strong than Theorem 
\ref{the5}. However, instead of containment (\ref{wat4}) it sets an inequality 
which is easy to check out.
\begin{corollary}\label{cor4}
Let $F_c$, $F_b$ and $F_a$ be the Fibonacci numbers where $c>b>a\geq 5$. Then a
numerical semigroup ${\sf S}\left(F_a,F_b,F_c\right)$ is symmetric if
\begin{eqnarray}
\lambda=\gcd(a,b)\geq 3\;,\;\;\;\gcd(\lambda,c)=1,2\;,\;\;\;F_cF_{\lambda}>
{\sf lcm}(F_a,F_b)-F_a-F_b\;.\label{wat9}
\end{eqnarray}
The Hilbert series and Frobenius number are given by (\ref{wat8}).
\end{corollary}
{\sf Proof} $\;\;\;$The two first relations in (\ref{wat9}) are taken from 
Theorem \ref{the5} and were proven in (\ref{wat6}). We have to use also the 
containment (\ref{wat4}). For this purpose take $F_c$ exceeding the Frobenius 
number of semigroup generated by two numbers $F_a/F_{\lambda}$ and $F_b/F_{
\lambda}$. This number ${\cal F}\left(F_a/F_{\lambda},F_b/F_{\lambda}\right)$ is
classically known due to Sylvester \cite{sy884}. So, we get
\begin{eqnarray}
F_c>\frac{F_a}{F_{\lambda}}\frac{F_b}{F_{\lambda}}-\frac{F_a}{F_{\lambda}}-
\frac{F_b}{F_{\lambda}}=\frac{{\sf lcm}(F_a,F_b)-F_a-F_b}{F_{\lambda}}\;,
\nonumber
\end{eqnarray}
where the Hilbert series $H\left(F_a,F_b,F_c\right)$ and Frobenius number ${\cal
F}\left(F_a,F_b,F_c\right)$ are given by (\ref{wat8}). Thus, Corollary is 
proven.$\;\;\;\;\;\;\Box$

We finish this Section by Example \ref{ex1} where the Fibonacci triple does 
satisfy the containment in (\ref{wat4}) but does not satisfy inequality in 
(\ref{wat9}).
\begin{example}$\{d_1,d_2,d_3\}=\{F_6=8,F_8=21,F_9=34\}$\label{ex1}
\begin{eqnarray}
&&\gcd(F_6,F_9)=F_3\;,\;\;\;\;\gcd(F_3,F_8)=1\;,\;\;\;\;F_8\in{\sf S}\left(
\frac{F_6}{F_3},\frac{F_9}{F_3}\right)={\sf S}\left(4,17\right)\;,\nonumber\\
&&f_1={\sf lcm}(F_6,F_9)=136\;,\;\;\;\;f_2=F_8\cdot F_3=42\;,\;\;\;\;
F_8\cdot F_3<{\sf lcm}(F_6,F_9)-F_6-F_9\;,\nonumber\\
&&H\left(F_6,F_8,F_9\right)=\frac{(1-z^{136})(1-z^{42})}{\left(1-z^8\right)
\left(1-z^{21}\right)\left(1-z^{34}\right)}\;,\;\;\;{\cal F}\left(F_6,F_8,F_9
\right)=115\;,\;\;\;\;G\left(F_6,F_8,F_9\right)=58\;.\nonumber
\end{eqnarray}
\end{example}
\section{Symmetric numerical semigroups generated by Lucas triple}\label{s5}
In this Section we consider symmetric numerical semigroups generated by three 
Lucas numbers $L_n$, $L_m$ and $L_k$, $n>m>k\geq 2$. Note that the case $k=2$ is
trivial because of Lucas number $L_2=3$ and Theorem \ref{the2}. The semigroup 
${\sf S}\left(L_2,L_m,L_n\right)$ is symmetric iff at least one of two 
requirements, $3\nmid L_m$ and $L_n\not\in{\sf S}\left(3,L_m\right)$, is broken.
\begin{theorem}\label{the6}
Let $L_k$, $L_m$ and $L_n$, $\;n,m,k\geq 3$, be the Lucas numbers and let 
\begin{eqnarray}
&&m=2^am^{\p}\;,\;\;n=2^bn^{\p}\;,\;\;k=2^ck^{\p}\;,\;\;\;\mbox{where}\;\;\;
m^{\p}=n^{\p}=k^{\p}=1\pmod 2\;,\;\;\;a,b,c\geq 0\;,\;\;\;\;\;\;\label{lis1a}\\
&&l=\gcd(m,n)=2^dl^{\p}\;,\;\;\;\mbox{where}\;\;\;l^{\p}=\gcd(m^{\p},n^{\p})=
1\pmod 2\;,\;\;\;d=\min\{a,b\}\;.\nonumber
\end{eqnarray}
Then a numerical semigroup generated by these numbers is symmetric iff $L_k$, 
$L_m$ and $L_n$ satisfy 
\begin{eqnarray}
L_k\in{\sf S}\left(\frac{L_m}{L_l},\frac{L_n}{L_l}\right)\;,\;\;\;\mbox{if}\;
\;a=b\;,\;\;\;\mbox{or}\;\;\;L_k\in{\sf S}\left(\frac{L_m}{2},\frac{L_n}{2}
\right)\;,\;\;\;\mbox{if}\;\;a\neq b\;,\label{lis2}
\end{eqnarray}
and one of three following relations:
\begin{eqnarray}
\begin{array}{ll}
1)\;\;\;\;a=b\neq 0\;,\;&\;\;a=b\neq c\;,\;\;\;3\nmid \gcd(k,l)\;,\\
2)\;\;\;\;a=b=0\;,\;\;\gcd\left(m^{\p},n^{\p}\right)>1\;\;,&\;
\left\{\begin{array}{ll}c=0\;,&\gcd\left(k^{\p},l^{\p}\right)=1\;,\\
c\neq 0\;,&3\nmid \gcd(k,l)\;,\end{array}\right.\\
3)\;\;\;\;a\neq b\;,\;\;\;3\mid\gcd(m,n)\;,\;&\;\;3\nmid k\;.\end{array}
\label{lis2a}
\end{eqnarray}
\end{theorem}
{\sf Proof} $\;\;\;$By Theorem \ref{the1} and Corollary \ref{cor1} a numerical  
semigroup ${\sf S}\left(L_k,L_m,L_n\right)$ is symmetric iff there exist 
two relatively not prime elements of its minimal generating set such that
\begin{eqnarray}
\eta=\gcd(L_n,L_m)>1\;,\;\;\gcd(L_k,\eta)=1\;,\;\;L_k\in{\sf S}\left(\frac{L_n}
{\eta},\frac{L_m}{\eta}\right)\;.\label{lis3}
\end{eqnarray}
Represent $n$ and $m$ as in (\ref{lis1a}) and substitute them into the 1st 
relation in (\ref{lis3}). By consequence of Theorem \ref{the4} it holds iff 
\begin{eqnarray}
1)\;\;a=b\;,\;\;\gcd(m,n)>1\;\;\;\;\;\mbox{or}\;\;\;\;\;2)\;\;a\neq 
b\;,\;\;3\mid\gcd(m,n)\;.\label{lis3a}
\end{eqnarray}
First, assume that the 1st requirement in (\ref{lis3a}) holds that results by 
Theorem \ref{the4} in $\eta=L_l$. Making use of notations (\ref{lis1a}) for $k$
move on to the 2nd requirement in (\ref{lis3}) and apply Corollary (\ref{cor2}).
Here we have to consider two cases $a=b\neq 0$ and $a=b=0$ separately.
\begin{eqnarray}
&&a=b\neq 0\;,\;\;\;a=b\neq c\;,\;\;\;3\nmid \gcd(k,l)=1\;,\label{lis3b}\\
&&a=b=0\;,\;\;\;\gcd\left(m^{\p},n^{\p}\right)>1\;,\;\;\;
\left\{\begin{array}{ll}c=0\;,&\gcd\left(k^{\p},l^{\p}\right)=1\;,\\
c\neq 0\;,&3\nmid\gcd(k,l)\;.\end{array}\right.\label{lis3c}
\end{eqnarray}
Now, assume that the 2nd requirement in (\ref{lis3a}) holds that results by
Theorem \ref{the4} in $\eta=2$. Making use of the 2nd requirement in 
(\ref{lis3}) and applying (\ref{wat3c}) we get,
\begin{eqnarray}
a\neq b\;,\;\;3\mid\gcd(m,n)\;,\;\;\;3\nmid k\;.\label{lis3d}
\end{eqnarray}
Combining (\ref{lis3b}), (\ref{lis3c}) and (\ref{lis3d}) we arrive at 
(\ref{lis2a}). The last requirement in (\ref{lis3}) together with Theorem 
\ref{the4} gives
\begin{eqnarray}
L_k=A\frac{L_m}{\eta}+B\frac{L_n}{\eta}=\left\{\begin{array}{lll}A\cdot L_m/L_l+
B\cdot L_n/L_l\;&\mbox{if}&\;a=b\\A\cdot L_m/2+B\cdot L_n/2\;&\mbox{if}&a\;
\neq b\end{array}\right.\;,\;\;\;A,B\in {\mathbb Z}_+\;,\nonumber
\end{eqnarray}
that proves (\ref{lis2}) and finishes proof of Theorem.$\;\;\;\;\;\;\Box$

By consequence of Theorem \ref{the6} the following Corollary holds for the 
most simple Lucas triples.
\begin{corollary}\label{cor5}
Let $L_{k^{\p}}$, $L_{m^{\p}}$ and $L_{n^{\p}}$ be the Lucas numbers with odd 
indices such that
\begin{eqnarray}
\gcd(m^{\p},n^{\p})>1\;,\;\;\;\;\gcd(m^{\p},n^{\p},k^{\p})=1\;.\label{lis5a}
\end{eqnarray}
Then a numerical semigroup generated by these numbers is symmetric iff 
\begin{eqnarray}
L_{k^{\p}}\in{\sf S}\left(\frac{L_{m^{\p}}}{L_{\gcd(m^{\p},n^{\p})}},
\frac{L_{n^{\p}}}{L_{\gcd(m^{\p},n^{\p})}}\right)\;.\label{lis5b}
\end{eqnarray}
\end{corollary}
Proof follows if we apply Theorem \ref{the6} in the case $a=b=c=0$, see 
(\ref{lis3c}).

We give without derivation the Hilbert series and Frobenius number for symmetric
semigroup ${\sf S}\left(L_{k^{\p}},L_{m^{\p}},L_{n^{\p}}\right)$.
\begin{eqnarray}
&&H\left(L_{n^{\p}},L_{m^{\p}},L_{k^{\p}}\right)=\frac{(1-z^{l_1})(1-z^{l_2})}
{\left(1-z^{L_{n^{\p}}}\right)\left(1-z^{m^{\p}}\right)\left(1-z^{L_{k^{\p}}}
\right)}\;,\;\;\;l_1=\frac{L_{n^{\p}}\cdot L_{m^{\p}}}{L_{\gcd(m^{\p},
n^{\p})}}\;,\nonumber\\
&&{\cal F}\left(L_{n^{\p}},L_{m^{\p}},L_{k^{\p}}\right)=l_1+l_2-(L_{n^{\p}}+
L_{m^{\p}}+L_{k^{\p}})\;,\;\;\;l_2=L_{k^{\p}}\cdot L_{\gcd(m^{\p},n^{\p})}
\;.\label{lis6}
\end{eqnarray}
In general, the containment (\ref{lis5b}) is hardly to verify because it 
presumes algorithmic procedure. Instead, one can formulate a simple inequality
which provide only the sufficient condition for semigroup ${\sf S}\left(L_{n^{
\p}},L_{m^{\p}},L_{k^{\p}}\right)$ to be symmetric.
\begin{corollary}\label{cor6}
Let $L_{n^{\p}}$, $L_{m^{\p}}$ and $L_{k^{\p}}$ be the Lucas numbers with odd 
indices such that (\ref{lis5a}) is satisfied and the following inequality holds,
\begin{eqnarray}
L_{k^{\p}}\;L_{\gcd(m^{\p},n^{\p})}>\frac{L_{n^{\p}}\;L_{m^{\p}}}
{L_{\gcd(m^{\p},n^{\p})}}-L_{n^{\p}}-L_{m^{\p}}\;.\label{lis7}
\end{eqnarray}
Then a numerical semigroup ${\sf S}\left(L_{n^{\p}},L_{m^{\p}},L_{k^{\p}}
\right)$ is symmetric and its Hilbert series and Frobenius number are 
given by (\ref{lis6}).
\end{corollary}
Its proof is completely similar to the proof of Corollary \ref{cor4} for 
symmetric semigroup generated by three Fibonacci numbers.

We finish this Section by Example \ref{ex2} where the Lucas triple does satisfy 
the containment in (\ref{lis5b}) but does not satisfy inequality (\ref{lis7}).
\begin{example}$\{d_1,d_2,d_3\}=\{L_{9}=76,L_{15}=1364,L_{17}=3571\}$\label{ex2}
\begin{eqnarray}
&&\gcd(L_9,L_{15})=L_3\;,\;\;\;\;\gcd(L_3,L_{17})=1\;,\;\;\;\;L_{17}\in{\sf S}
\left(\frac{L_9}{L_3},\frac{L_{15}}{L_3}\right)={\sf S}\left(19,341\right)\;,
\nonumber\\
&&l_1={\sf lcm}(L_9,L_{15})=25916\;,\;\;\;\;l_2=L_{17}\cdot L_3=14264\;,\;\;\;
\;L_{17}\cdot L_3<{\sf lcm}(L_9,L_{15})-L_9-L_{15}\;,\nonumber\\
&&H\left(L_9,L_{15},L_{17}\right)=\frac{(1-z^{25916})(1-z^{14264})}{\left(1-
z^{76}\right)\left(1-z^{1364}\right)\left(1-z^{3571}\right)}\;,\nonumber\\
&&{\cal F}\left(L_9,L_{15},L_{17}\right)=35189\;,\;\;\;\;
G\left(L_9,L_{15},L_{17}\right)=17595\;.\nonumber
\end{eqnarray}
\end{example}
\section*{Acknowledgement}
I thank C. Cooper for bringing the paper \cite{wmd91} to my attention.

\end{document}